\documentclass[a4paper,11pt]{amsart}

\usepackage{amsmath, amsthm, amssymb}
\usepackage{color}
\RequirePackage{textcmds}\relax
\ProvideTextCommandDefault{\cprime}{\tprime}
\usepackage[pagewise,mathlines]{lineno}
\usepackage[colorinlistoftodos,prependcaption,textsize=tiny]{todonotes}
\usepackage{cleveref}

\theoremstyle{plain}
\newtheorem{thm}{Theorem}[section]
\newtheorem{lem}[thm]{Lemma}
\newtheorem{prop}[thm]{Proposition}

\newtheorem{defi}{Definition}

\theoremstyle{remark}
\newtheorem{rem}[thm]{\bf{Remark}}

\begin{document}
\newcommand{\rr}{\mathbb{R}}
\newcommand{\nn}{\mathbb{N}}
\newcommand{\di}{\displaystyle}
  

\newcommand{\todoRL}[2][noinline]
{\todo[#1,linecolor=blue,backgroundcolor=green!20,bordercolor=black,caption={}]{\textbf{RL:}#2}}
\newcommand{\todoSF}[2][noinline]
{\todo[#1,linecolor=blue,backgroundcolor=yellow!20,bordercolor=black,caption={}]{\textbf{SF:}#2}}
\newcommand{\todoSA}[2][noinline]
{\todo[#1,linecolor=blue,backgroundcolor=red!20,bordercolor=black]{\textbf{SA:}#2}
}

\numberwithin{equation}{section}

                \title{Extension of the Best Polynomial Operator in Generalized Orlicz spaces}

                \author{ Sonia Acinas$^1$, Sergio Favier$^2$ and  Rosa Lorenzo$^3$ }
                       
			\address{$^1$
				 Departamento de Matem\'atica, Facultad de Ciencias Exactas Naturales, Universidad Nacional de La Pampa,
        L6300CLB Santa Rosa, La Pampa, Argentina}
        \email{sonia.acinas@gmail.com}

								\address{$^2$ 								
Instituto de Matem\'atica Aplicada San Luis, IMASL,
Universidad Nacional de San Luis and CONICET,
Ej\'ercito de los Andes 950,
D5700HHW San Luis, Argentina.
Departamento de Matem\'atica,
Universidad Nacional de San Luis, D5700HHW San Luis, Argentina.
        T.E. 54 266 4459649}
                \email{sfavier@unsl.edu.ar}
       
\address{$^3$ 								
Instituto de Matem\'atica Aplicada San Luis, IMASL,
Universidad Nacional de San Luis and CONICET,
Ej\'ercito de los Andes 950,
D5700HHW San Luis, Argentina.
Departamento de Matem\'atica,
Universidad Nacional de San Luis, D5700HHW San Luis, Argentina.
        T.E. 54 266 4459649}
                \email{rlorenzo@unsl.edu.ar}
				
				 \thanks{This paper was supported by Consejo Nacional de Investigaciones Científicas y Técnicas, Universidad Nacional de San Luis and Universidad Nacional de La Pampa. Grants: 11220200100694CO,   PROICO 3-0720 and PI 77(M).}

\begin{abstract}
In this paper, we consider the best multivalued polynomial approximation operator for functions in an Orlicz Space $L^{\varphi}(\Omega)$.
We obtain its characterization involving $\psi^-$ and $\psi^+$, which are the left and right derivative functions of $\varphi$. And then, we extend the operator to $L^{\psi^+}(\Omega)$.
We also get pointwise convergence of this extension,  where the Calderón-Zygmund class $t_m^p (x)$ adapted to  $L^{\psi^+}(\Omega)$ plays an important role. 
\end{abstract}

    \maketitle
    \markboth { \tiny{Extension of the Best Polynomial Approximation Operator in general Orlicz Spaces }}
        {Sonia Acinas, Sergio Favier and Rosa Lorenzo}
        \\



\section{Introduction}

Let $\mathcal{S}$ be the set of all non-decreasing functions $\psi:[0, \infty) \rightarrow[0, \infty)$ 
such that $\psi(0)=0$,  $\psi(s)>0$  if $s>0$,  $\psi(s) \rightarrow 0$ when $s \rightarrow 0^+$ and $\psi(s) \rightarrow \infty$ as $s \rightarrow \infty$.

We say that a function $\psi \in \mathcal{S}$ satisfies the $\Delta_2$ condition and we denote by $\psi \in \mathcal{S}\cap\Delta_2$, 
if there exists a constant $\Lambda=\Lambda_{\psi}>0$ such that
$\psi(2 s)\leqslant\Lambda_{\psi} \psi(s)$ \;for all $s\geq 0$.
Note that in this case, for any $C>0$ we have $ \psi ( C s) \leqslant K \psi (s),$ for a constant $K$ depending only on $C.$
In the rest of this work, for simplicity, we will write $K=\Lambda_{\psi}.$

Let $\Omega$ be a bounded Lebesgue measurable set in $\mathbb{R}^n$.  If $\psi \in
\mathcal{S}\cap\Delta_2$, we denote by $L^{\psi}(\Omega)$ the class of all Lebesgue real-valued measurable functions $f$ defined on $\Omega$ such that $\int_{\Omega}\psi(|f|)\,dx<\infty$, 
where $dx$ refers to Lebesgue measure in $\mathbb{R}^n$. In the sequel,
we sometimes write $\int_{\Omega}\psi(|f|)$ for $\int_{\Omega}\psi(|f|)\,dx$. 
And, if $A$ is a Lebesgue measurable set such that $A \subseteq \Omega$, we will denote by $|A|$ its Lebesgue measure 
and we will write $\int_A f$ for $\int_{\Omega \cap A} f$.

We set $\Phi$ for the class of all functions $\varphi$ 
defined on $[0,\infty)$ by $\varphi(s)=\int_{0}^{s} \psi(r) \,dr$, for some  $\psi \in \mathcal{S}$, 
which also satisfy the $\Delta_2$ condition.

Observe that if  $\varphi \in \Phi$, the space $L^{\varphi}(\Omega)$ is the classical Orlicz
space treated for example in \cite{KR61} and \cite{RR91}.

Let $\psi^{-}$, $\psi^{+}$ be the left and right derivatives of $\varphi$, respectively. We set  $\psi^{-}(0)=0$ and, due to $\frac{\varphi (s)}{s} \to 0,$ as $s\to 0^+$, we  have $\psi^{+}(0)=0.$

For $\varphi \in \Phi$, using the fact that $\varphi(s)=\int_{0}^{s} \psi^+(r) \,dr$  with
$\psi^+ \in \mathcal{S}$, it holds that
 \begin{equation}\label{eq:des-psi+ contra phi}
\frac{s}{2} \psi^{+}\left(\frac{s}{2}\right) \leqslant \varphi(s) \leqslant s\, \psi^{+}(s) \leqslant \varphi(2 s),\;\;
\mbox{ for all }\;\;s \geqslant 0.
\end{equation}
An analogous inequality is also satisfied by $\psi^-$ instead of $\psi^+$.
Now, from \eqref{eq:des-psi+ contra phi} and the $\Delta_2$ condition on $\varphi$, we get that $\psi^-$ and $\psi^+$ also satisfy the $\Delta_2$ condition.

Note that if $s<r$, due to the convexity of $\varphi$, we get  
\begin{equation}\label{eq:derivative-strictly-increasing}
\psi^+(s) \leqslant  \psi^-(r).
\end{equation}
Also, observe that by \eqref{eq:derivative-strictly-increasing}
and $\psi^-\in \Delta_2$, we have
\[
\psi^-(s) \leqslant  \psi^+(s) \leqslant  \Lambda_{\psi^-} \psi^-(s), \mbox{ for all  } s\geqslant 0;
\]
and, since $\psi^+$ is a non-decreasing function and $\psi^+\in \Delta_2$, we obtain
\begin{equation}\label{eq:subaditiva}
\frac{1}{2}(\psi^+(s)+\psi^+(r))\leqslant  \psi^+(s+r) \leqslant 
\Lambda_{\psi^+} (\psi^+(s)+\psi^+(r)),
\end{equation}
for every $s, r\geqslant 0.$
As before, \eqref{eq:subaditiva} holds for $\psi^-$ instead of $\psi^+$ and note that $L^{\varphi} (\Omega) \subseteq L^{\psi^-}(\Omega) = L^{\psi^+}(\Omega)$. 

Let $\Pi^m$ be the space of algebraic polynomials, defined on
$\mathbb{R}^n$, of degree at most $m.$ A polynomial  $P\in \Pi^m$ is said
to be  a best approximation of $f \in L^\varphi(\Omega)$   if and only if
\begin{equation}\label{defi-poli-aprox-bis}
\int_{\Omega} \varphi(|f-P|)\,dx =\inf_{Q\in \Pi^m} \int_{\Omega} \varphi(|f-Q|)\,dx.
\end{equation}
\begin{defi}
For $f\in L^{\varphi} (\Omega)$ we denote by $\mu_{\varphi}(f)$ the set of all
polynomials $P$ that satisfy  (\ref{defi-poli-aprox-bis}).
\end{defi}
Hereinafter we also refer to $\mu_{\varphi} (f)$ as the multivalued
operator defined for functions in $L^{\varphi} (\Omega)$ and images on
$\Pi^m.$

In \cite{AFZ14} the extension of the best polynomial
approximation operator $\mu_{\varphi}(f)$ is stated 
assuming that the function $\varphi$ is a differentiable
$N$-function, i.e. a differentiable convex function in $\Phi$ 
such that its derivative function $\varphi'$ satisfies
$\varphi'(0^+)=0$ and $\varphi'(s) \to \infty$ as $s \to \infty$.
The main goal of this paper
is to get a continuous extension and  properties of
$\mu_{\varphi}(f)$ when  $\varphi$ is not necessarily a differentiable $N$-function. 

We set this extension in Theorem \ref{thm:extension} on the assumption that $\varphi$  is not differentiable
and $\psi^-, \psi^+$ are the left and right derivatives of $\varphi \in  \Phi$ with  $\psi^-,\psi^+\in \Delta_2$. 
We point out that the techniques employed  to prove these results differ from those of \cite{AFZ14} due to the lack of continuity of the derivative of 
$\varphi$.

 We observe that the classical conditional expectation operator can be seen as the
 best approximation operator in the following setup. Let $f\in
 L^{2}$ defined on a probability space $(\Omega, {\mathcal A}, P)$ and
 let ${\mathcal A}_0 \subseteq {\mathcal A}$ be a subsigma algebra of ${\mathcal
 A}_0.$ Consider all ${\mathcal
 A}_0$ measurable functions which are also in $L^2.$ Then, we can see
 the conditional expectation as the best approximation operator
 defined on $L^2 (\Omega, {\mathcal A}, P)$ with images in $L^2 (\Omega, {\mathcal A}_0,
 P).$ 
 The monotone property allows us to continuously extend
 the operator to $L^1 (\Omega, {\mathcal A}, P),$ and this gives the classical
 conditional expectation.
 
 Further, for $\varphi (s) = s^p$ with $1<p<\infty,$ a similar extension to $L^{p-1}$ of the best approximation  operator has been considered in \cite{LR81}. Extension of
 this operator and another classical operator in  harmonic analysis
 have been also treated in general Orlicz spaces $L^{\varphi}$, see
 \cite{LR80}, \cite{BJ70} and \cite{KK91}.

 The special case of ${\mathcal A}_0 = \{\emptyset, \Omega\}$ takes us
 to the best approximation in $L^{\varphi} (\Omega)$ by constants
 functions and, then, to the extended operator, which has been extensively studied in
 \cite{MC01a}, \cite{FZ05} and \cite{FZ11}, where the differentiability of $\varphi$ is
 assumed. Also in \cite{CFZ08} the extension of the best approximation operator
 in $L^{\varphi} (\Omega)$ is treated for a general $\sigma-$lattice
 approximation class. Further, when the algebraic polynomials are
 considered as the approximation class, we refer to \cite{MC01a}, \cite{Cue11} and \cite{AFZ14},  where the best approximation operator $\mu_{\varphi}(f)$
 defined on $L^{\varphi} (\Omega)$ is extended to $L^{\varphi'}
 (\Omega),$ being $\varphi'$ the derivative function of
 $\varphi.$ Recently, in \cite{FGL22} an extension of this best approximation operator in
 Orlicz-Lorentz spaces was also considered and a generalization of these results can be found in  \cite{GKL23}.

 Subsequently, in \cite{FL2020}, taking the constant functions as the approximation class,
 it was obtained the extension of the best approximation operator, for the first time,  without the assumption of  the differentiability of $\varphi,$ and therefore the left and right
derivatives of $\varphi$ are involved.  The main objective of this manuscript is
 to obtain the extension of the best approximation operator  when $\varphi$ is an $N$-function without the differentiation hypothesis and where the approximation class is the algebraic polynomials.

In Section 2, we include an existence result of the best polynomial approximation in $L^{\varphi} (\Omega)$, to make this paper self-contained. In addition, we set a characterization of the best polynomial approximation operator for functions  $f \in L^{\varphi}(\Omega)$,   where $\varphi$ is not differentiable, and we derive some properties that will be useful to get the  extension of this operator for functions 
belonging to the bigger space $L^{\psi^{+}} (\Omega)$. 

In Section  3, we 
extend the best polynomial approximation operator from $L^{\varphi}(\Omega)$ to $L^{\psi^+}(\Omega)$, where we have to use a different method to that developed in \cite{AFZ14}, because
$\varphi$ is not a differentiable function in this article.

The uniqueness of the extended best polynomial approximation operator is studied in Section 4, and we also establish a sort of continuity of this operator.

Finally, in Section 5, we set pointwise convergence for the  extended best polynomial approximation  operator
for a wide class of functions $f$ closely related to the
Calder\'on-Zygmund class $t_m^p (x)$  introduced in 1961, see \cite{CZ61}. And, it becomes also a generalization of  \cite{CFZ12}.

\section{The best polynomial approximation operator in $L^{\varphi}(\Omega)$}

For $P\in \Pi^m$ we set $\|P\|_{\infty} = \sup\limits_{x\in \Omega}
|P(x)|$ and $\|P\|_1 = \int_{\Omega} |P|\, dx.$\\

First, we point out that the set $\mu_{\varphi}(f)$ is not empty for functions $f \in L^{\varphi}(\Omega)$, as stated in the following result. 

\begin{thm}
Let $\varphi \in \Phi$ 
 and let $f \in L^{\varphi}(\Omega)$. Then, there
exists $P\in \Pi^m$ such that
\[
\int_{\Omega} \varphi(|f-P|)\,dx =\inf_{Q\in \Pi^m} \int_{\Omega} \varphi(|f-Q|)\,dx.
\]
\end{thm}
The proof follows the same lines as Theorem 2.2 in \cite{AFZ14} and Theorem 2.1 in \cite{AF16}, where the differentiability of $\varphi$ is not used.

Theorem 2.3 of \cite{AFZ14} and Theorem 2.2 of \cite{AF16} give a characterization of the
best polynomial approximation for functions in $L^{\varphi}(\Omega)$ assuming that $\varphi\in\Phi$ has continuous derivative.  

Next, we obtain a characterization for the case of $\varphi$ not necessarily being a differentiable function.

\begin{thm}\label{thm:caract-bpa-non-differentiable}
Let $\varphi \in \Phi$ and let $f \in L^{\varphi}(\Omega)$.
Suppose that $\psi^-$ and $\psi^+$ are the left and right derivatives of $\varphi$, respectively.
\\ 
Then, $P\in
\Pi^m$ is in $\mu_{\varphi} (f)$ if and only if the following inequality holds
\begin{equation}\label{eq:caract-bpa-f>P}
\begin{aligned}
 \int_{\{f>P\} \cap\{Q>0\}} \psi^{-}(|f-P|) Q \,dx+\int_{\{f>P\} \cap\{Q<0\}} \psi^{+}(|f-P|) Q \,dx & \\
\leqslant \int_{\{f < P\} \cap\{Q<0\}} \psi^{-}(|f-P|) Q \,dx+\int_{\{f < P\} \cap\{Q>0\}} \psi^{+}(|f-P|) Q \,dx, &
\end{aligned}
\end{equation}
for all $Q\in \Pi^m.$ 
\end{thm}

\begin{proof}
For $R,S\in \Pi^m$, we define
\[
F_{R,S}(\epsilon)=\int_{\Omega} \varphi(|f-(R+\epsilon S)|) \,dx, 
\]
for $\epsilon\geqslant 0$.

Note that $F_{R,S}$ is a convex function on $[0, \infty)$. In fact,  for $a,b\geqslant 0$ such that $a+b=1$, we have
\[
\begin{aligned}
F_{R,S}\left(a \epsilon_{1}+b \epsilon_{2}\right) 
& \leqslant \int_{\Omega} \varphi\left(a\left|f-\left(R+\epsilon_{1} S\right)\right|+b\left|f-\left(R+\epsilon_{2} S\right)\right|\right) \,dx \\
& \leqslant a \int_{\Omega} \varphi\left(\left|f-\left(R+\epsilon_{1} S\right)\right|\right) dx+b \int_{\Omega} \varphi\left(\left|f-\left(R+\epsilon_{2} R\right)\right|\right) \,dx \\
&=a F_{R,S}\left(\epsilon_{1}\right)+b F_{R,S}\left(\epsilon_{2}\right),
\end{aligned}
\]
for every $\epsilon_{1}, \epsilon_{2} \geqslant 0$.

If $Q\equiv 0$, \eqref{eq:caract-bpa-f>P} holds trivially. So, throughout this proof, we consider the case $Q\not\equiv 0$.

Now, for $P \in \mu_{\varphi}(f)$  and  any $Q \in \Pi^{m}$,    it holds 
\[
F_{P,Q}(0)=
\min\limits _{\epsilon \in [0, \infty)} F_{P,Q}(\epsilon);
\]
and, this identity is satisfied if and only if $0 \leqslant F_{P,Q}^{\prime}\left(0^{+}\right)$, where  $F_{P,Q}^{\prime}\left(0^{+}\right)$ is the right derivative of $F_{P,Q}$ at  $\epsilon=0$. On the other hand, if for a fix $P\in \Pi^{m}$ it holds $F'_{P,Q} (0^{+}) \ge 0$ for every $Q\in \Pi^{m},$ then $P\in \mu_{\varphi}(f)$, where   
\[
\begin{aligned}
F_{P,Q}^{\prime}\left(0^{+}\right) &=\lim _{\epsilon \rightarrow 0^{+}} \frac{F_{P,Q}(0+\epsilon)-F_{P,Q}(0)}{\epsilon} \\
&=\lim _{\epsilon \rightarrow 0^{+}} \frac{1}{\epsilon}\left\{\int_{\Omega} \varphi(|f-(P+\epsilon Q)|) \,dx-\int_{\Omega} \varphi(|f-P|) \,dx\right\} .
\end{aligned}
\]
Let $F_1 =\{f>P+ \varepsilon Q\}$ and $F_2 =\{f\leqslant  P+ \varepsilon Q\}$. Taking into account that $F_1\cap \{Q\geqslant  0\}\cap \{f\leqslant  P\}= F_2\cap \{Q< 0\}\cap \{f\geqslant  P\} =\emptyset$ and the fact that Lebesgue measure of $\{Q=0\}$ is $0,$ we split $\Omega$ and obtain
\[
\begin{aligned}
F_{P,Q}^{\prime} & \left(0^{+}\right) =
\\
&\lim\limits_{\epsilon \rightarrow 0^{+}}
\left\{\int_{F_1  \cap\{Q>0\} \cap \{f>P\}}\left[\frac{\varphi(f-(P+\epsilon Q))-\varphi(f-P)}{\epsilon Q}\right] Q \,dx\right. 
\\
&+\int_{F_1 \cap\{Q<0\} \cap\{f>P\} }\left[\frac{\varphi(f-(P+\epsilon Q))-\varphi(f-P)}{\epsilon Q}\right] Q \,dx 
\\
&+\int_{ F_1 \cap\{Q<0\} \cap\{f \leqslant  P\} }\left[\frac{\varphi(f-(P+\epsilon Q))-\varphi(P-f)}{\epsilon Q}\right] Q \,dx 
\\
&+\int_{F_2 \cap\{Q>0\} \cap\{f <  P\} }\left[\frac{\varphi(P+\epsilon Q-f)-\varphi(P-f)}{\epsilon Q}\right] Q \,dx 
\\
&+\int_{F_2 \cap\{Q>0\} \cap\{f\geqslant P\} }\left[\frac{\varphi(P+\epsilon Q-f)-\varphi(f-P)}{\epsilon Q}\right] Q \,dx 
\\
&\left.+\int_{F_2 \cap\{Q<0\} \cap\{f < P\} }\left[\frac{\varphi(P+\epsilon Q-f)-\varphi(P-f)}{\epsilon Q}\right] Q \,dx\right\}
\\
&=\lim\limits_{\epsilon \rightarrow 0^{+}} 
\left\{I_1+I_2+I_3+I_4+I_5+I_6\right\}.
\end{aligned}
\]

For $ I_1,I_2, I_4$ and $I_6$, due to the convexity of $\varphi$ and (\ref{eq:subaditiva}), we have 
\begin{equation*}\label{acotacion}
\frac{|\varphi(|f-(P+\epsilon Q)|)- \varphi(|f-P|)|}{\epsilon }
\leqslant 
|Q| \Lambda_{\psi^+} \left(\psi^{+}(|f-P|)+\psi^{+}(|Q|)\right) ,
\end{equation*}
for  $0 < \epsilon \leqslant 1$. 
For $I_3$ and $I_5$, we get
\begin{equation}\label{acotacion-caso0}
\frac{|\varphi(|\epsilon Q)|)- \varphi(0)|}{\epsilon }
\leqslant   \\
|Q| \psi^{+}(|Q|),
\end{equation}
for $0 < \epsilon \leqslant 1$. 
And, observe that 
$|Q| \Lambda_{\psi^+} \left(\psi^{+}(|f-P|)+\psi^{+}(|Q|)\right) $
 is an integrable function on $\Omega$. 

Now, applying Dominated Convergence Theorem, we obtain 
\begin{equation}\label{eq:limI_1}
\lim\limits_{\epsilon \rightarrow 0^{+}} I_1=
-\int_{\{f > P\} \cap\{Q>0\}} \psi^{-}(|f-P|) Q \,dx,
\end{equation}
and 
\begin{equation}\label{eq:limI_2}
\lim\limits_{\epsilon \rightarrow 0^{+}} I_2=
-\int_{\{f > P\} \cap\{Q<0\}} \psi^{+}(|f-P|) Q \,dx.
\end{equation}

For any decreasing sequence $\left\{\epsilon_{n}\right\}$ that converges to 0, we consider the sets 
$A_{n}=\left\{f>P+\epsilon_{n} Q\right\} \cap$ $\{f<P\} \cap\{Q<0\}$. 
 Then  
$\left\{A_{n}\right\}$ is a decreasing sequence of sets such that, taking into account (\ref{acotacion-caso0}), we get  
\begin{equation}\label{eq:cota-integral-An}
\begin{aligned} \lim\limits_{\epsilon_n \rightarrow 0^{+}} \int_{A_{n}} \left[\frac{\varphi(f-(P+\epsilon_n Q))-\varphi(P-f)}{\epsilon_{n} Q}\right] &Q \,dx \leqslant  \\
&\lim\limits_{\epsilon_n \rightarrow 0^{+}} \int_{A_n}  \psi^{+}(|Q|)|Q|\, dx.
\end{aligned}
\end{equation}
For any measurable set $A\subseteq \Omega$, we set 
$\mu(A)=\int_A \psi^+(|Q|)|Q|\,dx$.
Then 
$\mu(A_n)<\infty$ and due to $\bigcap\limits_{n=1}^{\infty} A_{n}=\emptyset$, we get
$\lim\limits_{n \rightarrow \infty} \mu\left(A_{n}\right) = 
 0$, and therefore
\begin{equation}\label{eq:limI_3}
\lim\limits_{\epsilon \rightarrow 0^{+}} I_3=0.
\end{equation}

Again, by Dominated Convergence Theorem and $\psi^+(0)=0,$ we have
\begin{equation}\label{eq:limI_4}
\lim\limits_{\epsilon \rightarrow 0^{+}} I_4=
\int_{\{f < P\} \cap\{Q>0\}} \psi^{+}(|f-P|) Q \,dx,
\end{equation}
and
\begin{equation}\label{eq:limI_6}
\lim\limits_{\epsilon \rightarrow 0^{+}} I_6=
\int_{\{f < P\} \cap\{Q<0\}} \psi^{-}(|f-P|) Q \,dx.
\end{equation}

Using a similar argument as for reaching \eqref{eq:cota-integral-An}, we consider the sets  
$B_{n}=\left\{f \leqslant P+\epsilon_{n} Q\right\} \cap\{f>P\} \cap\{Q>0\}$, 
where
$\bigcap\limits_{n=1}^{\infty}\left\{f \leqslant P+\epsilon_{n} Q\right\} \cap\{f>P\} \cap\{Q>0\}=\emptyset$,  
and due to $\psi^{+} (0) =0$, we get  
\begin{equation*}\label{partedelim5}
\lim\limits_{\epsilon \rightarrow 0^{+}} 
\int_{\{f \leqslant P + \varepsilon_n Q\} \cap \{f=P\} \cap \{Q>0\}} \frac{\varphi(\varepsilon  Q)}{\varepsilon Q} Q \,dx =0.
\end{equation*}

Thus, we obtain 
\[
\lim\limits_{\epsilon \rightarrow 0^{+}} I_5=0.
\]

From \eqref{eq:limI_1} to \eqref{eq:limI_6}, 
and taking into account that $F_{P,Q}^{\prime}\left(0^{+}\right) \geqslant 0$, 
we get
\[
\begin{aligned}
& \int_{ \{f>P\} \cap\{Q>0\}} \psi^{-}(|f-P|) Q \,dx
+\int_{\{f>P\} \cap\{Q<0\}} \psi^{+}(|f-P|) Q \,dx \\
\leqslant & \int_{\{f < P\} \cap\{Q<0\}} \psi^{-}(|f-P|) Q \,dx
+\int_{ \{f < P\} \cap\{Q>0\}} \psi^{+}(|f-P|) Q \,dx,
\end{aligned}
\]
which is inequality  \eqref{eq:caract-bpa-f>P}.
\end{proof}

\begin{rem}\label{Rem1}
Taking $-Q$ instead of $Q$  in inequality  \eqref{eq:caract-bpa-f>P}, we get 
\begin{equation}\label{eq:caract-bpa-f<P}
\begin{aligned}
\int_{\{f<P\} \cap\{Q>0\}} \psi^{-}(|f-P|) Q \,dx+\int_{\{f< P\} \cap\{Q<0\}} \psi^{+}(|f-P|) Q \,dx &\\
\leqslant \int_{\{f > P\} \cap\{Q<0\}} \psi^{-}(|f-P|) Q \,dx+\int_{\{f > P\} \cap\{Q>0\}} \psi^{+}(|f-P|) Q \,dx,&
\end{aligned}
\end{equation}
for any $Q \in \Pi^m$. Now, 
writing \eqref{eq:caract-bpa-f<P} in terms of $|Q|$, 
we obtain
\begin{equation}\label{remakb}
\begin{aligned}
\int_{\{f > P\} \cap\{Q<0\}} \psi^{-}(|f-P|)|Q| \,dx+\int_{\{f<P\} \cap\{Q>0\}} \psi^{-}(|f-P|)|Q|\, dx& \\
\leqslant \int_{\{f< P\} \cap\{Q<0\}} \psi^{+}(|f-P|)|Q| \,dx+\int_{\{f > P\} \cap\{Q>0\}} \psi^{+}(|f-P|)|Q|\,dx,&
\end{aligned}
\end{equation}
for all $Q\in \Pi^m,$ and taking $-Q$ instead of $Q$ in \eqref{remakb}, we also get
\begin{equation}\label{Remark}
\begin{aligned}
\int_{\{f > P\} \cap\{Q>0\}} \psi^{-}(|f-P|)|Q| \,dx+\int_{\{f<P\} \cap\{Q<0\}} \psi^{-}(|f-P|)|Q|\, dx& \\
\leqslant \int_{\{f< P\} \cap\{Q>0\}} \psi^{+}(|f-P|)|Q| \,dx+\int_{\{f > P\} \cap\{Q<0\}} \psi^{+}(|f-P|)|Q|\,dx,&
\end{aligned}
\end{equation}
for all $Q\in \Pi^m.$
\end{rem}

The following result, in the spirit of Theorem 2.1 of \cite{CFZ12}, Theorem 2.4 of \cite{AFZ14} and
Theorem 2.4 of \cite{AF16},
gives us an inequality that will become a useful tool for the extension of the best approximation operator.

\begin{thm}\label{thm:des-map-phi-non differentiable}
Let $\varphi \in \Phi$  and let $f \in L^{\varphi}(\Omega).$
Suppose that the polynomial $P\in \Pi^m$ is  in $\mu_{\varphi} (f).$
Then
\begin{equation}\label{eq:control-mape-by-f}
\int_{\Omega} \psi^{-}(|P|)|P| \,dx 
\leqslant 5 \Lambda_{\psi^-}\|P\|_{\infty} \int_{\Omega} \psi^{+}(|f|) \,dx.
\end{equation}
\end{thm}

\begin{proof}
Due to the monotonicity of the function $\psi^{-}$, we have
\begin{equation}\label{eq:subaditiva-psi-menos}
\int_{\Omega} \psi^{-}(|P|)|P| \,dx \leqslant \int_{\Omega} \psi^{-}(|f-P|+|f|)|P| \,dx.
\end{equation}
As $\psi^{-} \in \Delta_2$, from \eqref{eq:subaditiva}  there  exists $\Lambda_{\psi^{-}}>0$ such that 
\begin{equation}\label{eq:|f|P por Delta 2}
\int_{\Omega} \psi^{-}(|f-P|+|f|)|P| \,dx \leqslant \Lambda_{\psi^{-}} \int_{\Omega} \psi^{-}(|f-P|)|P| d x+ \Lambda_{\psi^{-}} \int_{\Omega} \psi^{-}(|f|)|P| \,dx.
\end{equation}
We  write
\begin{equation}\label{eq:|f-P|P sobre Omega}
\begin{aligned}
&\int_{\Omega} \psi^{-}(|f-P|)|P| \,dx 
\\
=&\int_{\{f>P\} \cap\{P>0\}} \psi^{-}(|f-P|) P dx+\int_{\{f>P\} \cap\{P<0\}} \psi^{-}(|f-P|)(-P) dx \\
+&\int_{\{f<P\} \cap\{P>0\}} \psi^{-}(|f-P|) P dx+\int_{\{f<P\} \cap\{P<0\}} \psi^{-}(|f-P|)(-P) dx.
\end{aligned}
\end{equation} 
We have
\begin{equation}\label{eq:|f-P|P sobre f>P}
\int_{\{f>P\} \cap\{P>0\}} \psi^{-}(|f-P|) P \,dx \leqslant \int_{\{f>P\} \cap\{P>0\}} \psi^{-}(|f|) P \,dx, 
\end{equation}
and
\begin{equation}\label{eq:|f-P|P sobre f<P}
\int_{\{f<P\} \cap\{P<0\}} \psi^{-}(|f-P|)(-P) \,dx \leqslant \int_{\{f<P\} \cap\{P<0\}} \psi^{-}(|f|)(-P) \,dx. 
\end{equation}
From inequality \eqref{remakb} with $Q=P$, we obtain 
\begin{equation}\label{eq:|f-P|(-P) sobre f<P}
\begin{aligned}
&\int_{\{f > P\} \cap\{P<0\}} \psi^{-}(|f-P|)(-P) dx
+\int_{\{f<P\} \cap\{P>0\}} \psi^{-}(|f-P|) P dx \\
\leqslant &\int_{\{f < P\} \cap\{P<0\}} \psi^{+}(|f-P|)(-P) dx+\int_{\{f > P\} \cap\{P>0\}} \psi^{+}(|f-P|) P dx \\
\leqslant &\int_{ \{f<P\} \cap\{P<0\}} \psi^{+}(|f|)(-P) dx+\int_{\{f > P\} \cap\{P>0\}} \psi^{+}(|f|) P dx.
\end{aligned}
\end{equation}
Then, by \eqref{eq:|f-P|P sobre Omega}, \eqref{eq:|f-P|P sobre f>P}, \eqref{eq:|f-P|P sobre f<P} and \eqref{eq:|f-P|(-P) sobre f<P}, we get 
\begin{equation}\label{eq: de |f|P sobre Omega}
\begin{aligned}
&\int_{\Omega} \psi^{-}(|f-P|)|P| \,dx \\
 \leqslant &\int_{\{f>P\} \cap\{P>0\}} \psi^{-}(|f|) P \,dx+\int_{\{f<P\} \cap\{P<0\}} \psi^{-}(|f|)(-P) \,dx \\
+&\int_{\{f > P\} \cap\{P>0\}} \psi^{+}(|f|) P \,dx+\int_{\{f<P\} \cap\{P<0\}} \psi^{+}(|f|)(-P) \,dx .
\end{aligned}
\end{equation}
From \eqref{eq:subaditiva-psi-menos}, \eqref{eq:|f|P por Delta 2}, \eqref{eq: de |f|P sobre Omega} and the fact that  $\psi^{-} \leqslant \psi^{+}$, 
we obtain 
\[
\begin{aligned}
&\int_{\Omega} \psi^{-}(|P|)|P| \,dx 
\\
 \leqslant 2 \Lambda_{\psi^{-}} &\int_{\{f>P\} \cap\{P>0\}} \psi^{+}(|f|) P \,dx
+2 \Lambda_{\psi^{-}} \int_{\{f<P\} \cap\{P<0\}} \psi^{+}(|f|)(-P) \,dx \\
+\Lambda_{\psi^{-}} &\int_{\Omega} \psi^{+}(|f|)|P| \,dx 
 \leqslant 5 \Lambda_{\psi^{-}}\int_{\Omega} \psi^{+}(|f|)|P| \,dx.
\end{aligned}
\]
Then 
\[
\int_{\Omega} \psi^{-}(|P|)|P| \,dx \leqslant 5 \Lambda_{\psi^{-}} \int_{\Omega} \psi^{+}(|f|)|P| \,dx \leqslant 5 \Lambda_{\psi^{-}}\|P\|_{\infty} \int_{\Omega} \psi^{+}(|f|) \,dx. 
\]
\end{proof}

\begin{rem}\label{rem:control-mape-by-f-Lpsi}
To prove Theorem \ref{thm:des-map-phi-non differentiable}, we have used that the polynomial $P$ is a solution of inequality \eqref{eq:caract-bpa-f>P} and this inequality makes sense for $f \in L^{\psi^+}(\Omega)$.
Therefore, inequality \eqref{eq:control-mape-by-f} holds for any polynomial $P$ that satisfies \eqref{eq:caract-bpa-f>P}, even  if $f \in L^{\psi^+}(\Omega)$.
\end{rem}

\section{Extension of the best polynomial approximation operator}

In Theorem 3.3 of \cite{AFZ14} and  Theorem 3.3 of \cite{AF16} it was proved that if $\varphi$ is differentiable, then $\mu_{\varphi'}(f)$ is a non-empty set;  and, it was established  
the extension of the
operator $\mu_{\varphi}(f)$ from the space $L^{\varphi}(\Omega)$ to a bigger
space $L^{\varphi'}(\Omega),$ for an $N$-function $\varphi$ given by
$\varphi(s)=\int_0^s \varphi'(r) \, dr$. 
In this section, we  will get the  extension of $\mu_{\varphi} (f)$ 
without requesting the function $\varphi$ to be differentiable.
That is, we will obtain the extension of $\mu_{\varphi}(f)$ 
from the space $L^{\varphi}(\Omega)$ into the bigger space $L^{\psi^+}(\Omega)$.

\begin{defi}\label{defi:extended-bpa}
For  $\varphi \in \Phi$ and for a function $f$  in   $ L^{\psi^{+}}(\Omega),$ we say that a polynomial $P$ is a best extended approximation polynomial, 
if $P$ is a solution of the next inequality 
\begin{equation}\label{eq:def-extended-q-neg}
\begin{aligned}
&\int_{\{f>P\} \cap\{Q>0\}} \psi^{-}(|f-P|)|Q| \,dx+\int_{\{f < P\} \cap \{Q<0\}} \psi^{-}(|f-P|)|Q| \,dx \\
\leqslant &\int_{\{f < P\} \cap\{Q>0\}} \psi^{+}(|f-P|)|Q| \,dx+\int_{\{f>P\} \cap\{Q<0\}} \psi^{+}(|f-P|)|Q| \,dx,
\end{aligned}
\end{equation}
for all $Q\in \Pi^m.$ 
\end{defi}

\begin{defi}\label{def:extended-bphi-operator}
For $f \in L^{\psi^{+}}(\Omega)$, we denote by $\mu_{\psi^{+}}^{\Omega} (f)= \mu_{\psi^{+}}(f)$ the set of polynomials that satisfies \eqref{eq:def-extended-q-neg}.
The application $\mu_{\psi^{+}}:  L^{\psi^{+}} \to 2^{\Pi^{m}}$ is called the extended best $\varphi-$approximation polynomial operator.
\end{defi}

\begin{rem}\label{positive-measure}
We point out that if $f \in L^{\psi^{+}}(\Omega)$ and $f$ does not belong to  $\Pi^{m},$  then for $P\in \mu_{\psi^{+}}(f)$ and using \eqref{eq:def-extended-q-neg}, we deduce that   $\{f>P\}$ and $\{f<P\}$ have positive measures.
In fact, if, for example, $\{f>P\}$ has zero measure, choosing $Q\equiv -1$ in \eqref{eq:def-extended-q-neg}, it leads to a contradiction. The case  $|\{f<P\}|=0$ follows similarly with $Q\equiv 1$.
\end{rem}

Before we treat the existence of the extended operator, we set some auxiliary results that will be needed later.

\begin{lem}\label{lem:1st-auxiliar-caract-extended}
Let $\varphi \in \Phi$ and suppose that  $\left\{f_{n}\right\}$ is a sequence in  $L^{\varphi}(\Omega)$ such that there exists a positive constant  $C$ 
satisfying $\int_{\Omega} \psi^{+}\left(\left|f_{n}\right|\right) \,dx \leqslant C$. 
Then, the set $\left\{\|P\|_{\infty}\right.$ : $\left.P \in \mu_{\varphi}\left(f_{n}\right), n=1,2, \ldots\right\}$ is bounded.
\end{lem}

The proof of Lemma \ref{lem:1st-auxiliar-caract-extended} follows the same lines as Lemma 3.1 of  \cite{AFZ14}, where the differentiability of $\varphi$ was not used.

\begin{lem}\label{eq:conv in L1 of  psi+fn}
Let $\varphi \in \Phi$, suppose that $f_{n}$ and $f$ are functions in $L^{\psi^{+}}(\Omega)$ such that
\[
\int_{\Omega} \psi^{+}\left(\left|f_{n}-f\right|\right) \,dx \rightarrow 0,
\]
as $n \rightarrow \infty$, and  assume that there exists $H \in L^{\psi^+}(\Omega)$ such that $|f_n (x)|\leqslant H(x)$ a.e. for $x \in \Omega$.

Besides, suppose that $g_{n}$ and $g$ are measurable functions such that $\left|g_{n}\right| \leqslant C$ for all $n \in \mathbb{N}$
and $g_{n} \rightarrow g$ a.e. for $x \in \Omega \cap\{f \neq 0\}$. Then, there exists a  subsequence $\left\{n_{k}\right\}$ such that
\[
\begin{split}
 \int_{\Omega} \psi^{-}(|f|) g \,dx &\leqslant
 \liminf \int_{\Omega} \psi^{-}\left(\left|f_{n_{k}}\right|\right) g_{n_{k}} \,dx  
 \\
 &\leqslant
 \limsup 
 \int_{\Omega} \psi^{+}\left(\left|f_{n_{k}}\right|\right) g_{n_{k}} \,dx 
 \leqslant  
 \int_{\Omega} \psi^{+}(|f|) g \,dx,
 \end{split}
 \]
as $k \rightarrow \infty$.
\end{lem}

\begin{proof}
Since $\psi^{+}$ is a non-decreasing function and $\psi^{+} (s) >0,$ for $s>0,$ there exists  a subsequence $\{f_{n_k}\}$ which converges to $f$ a.e. on $\Omega$. 
Next, we use  Fatou's Lemma and lower semicontinuity of $\psi^-$ for the first inequality. 
And, the third inequality follows using the facts that $|f_n(x)|\leq H(x)$ a.e. $x\in \Omega$ and $\left|g_{n}\right| \leqslant C$ for all $n \in \mathbb{N}$, the upper semicontinuity of 
$\psi^+$ and  Fatou Lebesgue's Theorem.
\end{proof}

Now, we  extend the best polynomial approximation operator from  $L^{\varphi}$ into the bigger space  $L^{\psi^{+}}.$

\begin{thm}\label{thm:extension}
Let $\varphi \in \Phi$  and let  $f \in L^{\psi^{+}}(\Omega)$. 
Then, there exists a polynomial $P,$ which satisfies inequality
\eqref{remakb} of Remark \ref{Rem1}.

In addition, it holds 
\begin{equation}\label{eq:inequality-extended}
\int_{\Omega} \varphi(|P|) \,dx \leqslant 5\Lambda_{\psi^+}\|P\|_{\infty} \int_{\Omega} \psi^{+}(|f|) \,dx.
\end{equation}
\end{thm}

\begin{proof}
Given $f   \in L^{\psi^{+}}(\Omega)$, we define $f_n=\min(\max(f,n),-n)$ with  $n \in \mathbb{N}$, 
then  $-n \leqslant f_{n} \leqslant n$ on $\Omega$, $f_{n} \in L^{\varphi}(\Omega)$ for every $n \in \mathbb{N}$, 
$\lim\limits_{n \rightarrow \infty} f_{n}=f$ 
and 
\begin{equation}\label{eq:fn-bounded-by-f}
\left|f_{n}\right| \leqslant|\max(f,-n)| \leqslant|f|.
\end{equation}
Then, by Theorem \ref{thm:caract-bpa-non-differentiable}, there exists $P_n \in \mu_{\varphi}(f_n)$, for every $n \in \mathbb{N}$. 
As $f \in L^{\psi^{+}}(\Omega)$ and \eqref{eq:fn-bounded-by-f} holds, there exists $C>0$ such that 
$\int_{\Omega} \psi^{+}\left(\left|f_{n}\right|\right) \leqslant \int_{\Omega} \psi^{+}(|f|) \leqslant C$, for every $n \in \mathbb{N}$. 
Now, by Lemma \ref{lem:1st-auxiliar-caract-extended}, the sequence $\left\{P_{n}\right\}_{n \in \mathbb{N}}$ is uniformly bounded. 
Then, there exists a subsequence  of polynomials $\left\{P_{n}\right\}$  uniformly convergent to $P\in \Pi^{m}$ over $\Omega$. 

For simplicity, we call $\left\{P_{n}\right\}$ such a subsequence  and it holds $P_{n} \rightarrow P$ when $n \rightarrow \infty$. 

For $g \in L^{\psi^{+}}(\Omega)$ and $S\in \Pi^m$, we set 
$F^{-} (g,S,Q) = \psi^{-} (|g-S|) |Q|$  and $F^{+} (g,S,Q) = \psi^{+} (|g-S|) |Q|$.
As $P_{n} \in \mu_{\varphi}\left(f_{n}\right)$, by Theorem \ref{thm:caract-bpa-non-differentiable} and inequality \eqref{remakb} of Remark \ref{Rem1}, we have
\begin{equation}\label{eq:consecuencia-caract-bpa}
\begin{aligned}
\int_{\left\{f_{n} > P_{n}\right\} \cap\{Q<0\}} F^{-} (f_n,P_n,Q) 
&+\int_{\left\{f_{n}<P_{n}\right\} \cap\{Q>0\}} F^{-} (f_n,P_n,Q) \\
\leqslant 
\int_{\left\{f_{n} > P_{n}\right\} \cap\{Q>0\}} F^{+} (f_n,P_n,Q)
&+
\int_{\left\{f_{n} < P_{n}\right\} \cap\{Q<0\}} F^{+} (f_n,P_n,Q),
\end{aligned}
\end{equation}
for every $Q \in \Pi^{m}$.

We set
$l=\liminf\limits_{n \rightarrow \infty}$ and  $L= \limsup\limits_{n \rightarrow \infty}.$
Now, we use \eqref{eq:fn-bounded-by-f} and Lemma \ref{lem:1st-auxiliar-caract-extended} to get subsequences $f_{n_k}, P_{n_k}$ written again by $f_n$, $P_n$ such that \\
$\int_{\Omega} \psi^{+}\left(\left|(f_{n}-P_n)-(f-P)\right|\right) \,dx \rightarrow 0,$
and also  $\chi_{\{f_n> P_n\}}$ and $\chi_{\{f_n < P_n\}}$
converge to $\chi_{\{f> P\}}$ and $\chi_{\{f < P\}}$, 
respectively, for $x \in \{f\neq P\}$. Then, 
  Lemma \ref{eq:conv in L1 of  psi+fn} gives
\[
\begin{aligned}
&
\int_{\{f > P\} \cap\{Q<0\}} F^{-} (f,P,Q)  
+\int_{\{f<P\} \cap\{Q>0\}} F^{-} (f,P,Q) 
\\
\leqslant l &\int_{\left\{f_{n}>P_{n}\right\} \cap\{Q<0\}} F^{-} (f_n,P_n,Q)  
+ l
\int_{\left\{f_{n}<P_{n}\right\} \cap\{Q>0\}} F^{-} (f_n,P_n,Q) 
\\
\leqslant l &
\left( \int_{\{f_{n}>P_{n}\} \cap\{ Q<0 \}} F^{-} (f_n,P_n,Q)   
+ \int_{\{f_{n}<P_{n}\} \cap\{Q>0\}} F^{-} (f_n,P_n,Q)\right) 
\\
\leqslant l & \left(\int_{\left\{f_{n}>P_{n}\right\} \cap\{Q>0\}} F^{+} (f_n,P_n,Q)   + \int_{\left\{f_{n} < P_{n}\right\} \cap\{Q<0\}} F^{+} (f_n,P_n,Q)   \right)
\\
\leqslant L &
\left(\int_{\left\{f_{n}>P_{n}\right\} \cap\{Q>0\}} F^{+} (f_n,P_n,Q)  
+ \int_{\left\{f_{n} < P_{n}\right\} \cap\{Q<0\}} F^{+} (f_n,P_n,Q) \right)
\\
\leqslant L& \int_{\left\{f_{n}>P_{n}\right\} \cap\{Q>0\}} F^{+} (f_n,P_n,Q) 
+ L 
\int_{\left\{f_{n} < P_{n}\right\} \cap\{Q<0\}} F^{+} (f_n,P_n,Q)
\\
\leqslant
& \int_{\left\{f>P\right\} \cap\{Q>0\}} F^{+} (f,P,Q) 
+  
\int_{\left\{f < P\right\} \cap\{Q<0\}} F^{+} (f,P,Q), 
\end{aligned}\]
for every $Q \in \Pi^{m}$.
Thus, we have obtained inequality \eqref{remakb} of Remark \ref{Rem1}.

Finally, by Remark \ref{rem:control-mape-by-f-Lpsi} and Theorem \ref{thm:des-map-phi-non differentiable}, we get
\[
\int_{\Omega} \varphi(|P|) \,dx 
\leqslant \int_{\Omega} |P|\psi^+(|P|) \,dx
\leqslant 5 \Lambda_{\psi^+} \|P\|_{\infty} \int_{\Omega} \psi^{+}(|f|) \,dx.
\]
\end{proof}

\section{Properties of $\mu_{\psi^{+}}(f)$}

In this section, we set some properties of the extended best polynomial approximation operator.
These results establish a remarkable difference with the case of constants as the approximation class considered in \cite{FL2020}.

In Theorem \ref{thm:extension}, we have proved the existence of the extended best polynomial approximation operator. 
Next, we provide a result stating the conditions under which the uniqueness is achieved.

\begin{thm}\label{thm:unique extended bpao}
Let  $\varphi \in \Phi$  and let $f \in L^{\psi^{+}}(\Omega)$. 
\\
If  $\psi^{+}$ is a strictly increasing function, then
there exists a unique extended best polynomial  approximation for each $f \in L^{\psi^{+}}(\Omega)$. 
\end{thm}

\begin{proof}
If $f \in \Pi^m$, the result follows straightforwardly from  \eqref{defi-poli-aprox-bis}.  

Now, from Remark \ref {positive-measure},  for any $P \in \mu_{\psi^{+}}(f)$ we have  $|\{f>P\}| >0.$ Let $P_1,P_2 \in \mu_{\psi^{+}}(f)$ and assume that $P_1$ and $P_2$ are different. 
Then, $Q=P_2-P_1$ is not the zero polynomial. 

First, we  deal  with the case
$\left|\left\{f>P_1\right\} \cap\left\{P_2>P_1\right\}\right| >0 $. 
\\
From \eqref{eq:def-extended-q-neg}, the strict monotonicity of $\psi^{-}$ and  $\psi^+$, \eqref{eq:derivative-strictly-increasing},  and \eqref{remakb},
 we get 
\[
\begin{split}
&\int_{\left\{f>P_1\right\} \cap\left\{P_2>P_1\right\}} \psi^{-}\left(\left|f-P_1\right|\right)|Q|\,dx+\int_{\left\{f<P_1\right\} \cap\left\{P_1>P_2\right\}} \psi^{-}\left(\left|f-P_1\right|\right)|Q|\,dx \leqslant  \\
&\int_{\left\{f< P_1\right\} \cap\left\{P_2>P_1\right\}} \psi^{+}\left(\left|f-P_1\right|\right)|Q|\,dx+\int_{\left\{f>P_1\right\} \cap\left\{P_1>P_2\right\}} \psi^{+}\left(\left|f-P_1\right|\right)|Q|\,dx<\\
&\int_{\left\{f<P_2\right\} \cap\left\{P_2>P_1\right\}} \psi^{-}\left(\left|f-P_2\right|\right)|Q|\,dx+\int_{\left\{f>P_2\right\} \cap\left\{P_1>P_2\right\}} \psi^{-}\left(\left|f-P_2\right|\right)|Q|\,dx\leqslant  
\\
&\int_{\left\{f>P_2\right\} \cap\left\{P_2>P_1\right\}} \psi^{+}\left(\left|f-P_2\right|\right)|Q|\,dx+\int_{\left\{f < P_2\right\} \cap\left\{P_1>P_2\right\}} \psi^{+}\left(\left|f-P_2\right|\right)|Q|\,dx<
\\
&\int_{\left\{f>P_1\right\} \cap\left\{P_2>P_1\right\}} \psi^{-}\left(\left|f-P_1\right|\right)|Q|\,dx+\int_{\left\{f<P_1\right\} \cap\left\{P_1>P_2\right\}} \psi^{-}\left(\left|f-P_1\right|\right)|Q|\,dx,
\end{split}\]
which is a contradiction.

Now, suppose  
$|\{f>$ $\left.P_1\right\} \cap\left\{P_2<P_1\right\}| >0$.
\\
Again, from the strict monotonicity of $\psi^{-}$  and $\psi^+$, \eqref{eq:derivative-strictly-increasing},  
\eqref{remakb} and \eqref{eq:def-extended-q-neg}, we have
\[
\begin{split}
&\int_{\left\{f>P_2\right\} \cap\left\{P_2<P_1\right\}} \psi^{-}\left(\left|f-P_2\right|\right)|Q|\,dx +\int_{\left\{f<P_2\right\} \cap\left\{P_2>P_1\right\}} \psi^{-}\left(\left|f-P_2\right|\right)|Q|\,dx  \leqslant 
\\
&\int_{\left\{f< P_2\right\} \cap\left\{P_2<P_1\right\}} \psi^{+}\left(\left|f-P_2\right|\right)|Q|\,dx +\int_{\left\{f>P_2\right\} \cap\left\{P_2>P_1\right\}} \psi^{+}\left(\left|f-P_2\right|\right)|Q|\,dx <
\\
&\int_{\left\{f<P_1\right\} \cap\left\{P_2<P_1\right\}} \psi^{-}\left(\left|f-P_1\right|\right)|Q|\,dx +\int_{\left\{f>P_1\right\} \cap\left\{P_2>P_1\right\}} \psi^{-}\left(\left|f-P_1\right|\right)|Q|\,dx  \leqslant 
\\
&\int_{\left\{f>P_1\right\} \cap\left\{P_2<P_1\right\}} \psi^{+}\left(\left|f-P_1\right|\right)|Q|\,dx +\int_{\left\{f < P_1\right\} \cap\left\{P_2>P_1\right\}} \psi^{+}\left(\left|f-P_1\right|\right)|Q|\,dx <
\\
&\int_{\left\{f>P_2\right\} \cap\left\{P_2<P_1\right\}} \psi^{-}\left(\left|f-P_2\right|\right)|Q|\,dx 
+\int_{\left\{f<P_2\right\} \cap\left\{P_2>P_1\right\}} \psi^{-}\left(\left|f-P_2\right|\right)|Q|\,dx ,
\end{split}
\]
which is also a contradiction.
\end{proof}

\begin{thm}\label{thm:prop-traslacion}
Let $\varphi \in \Phi.$ 
If  $f \in L^{\psi^{+}}(\Omega)$, then $\mu_{\psi^{+}}(f+P)=$ $\mu_{\psi^{+}}(f)+P$ for every $P \in \Pi^{m}$.
\end{thm}

\begin{proof}
For $f+P \in L^{\psi^{+}}(\Omega)$ and $\tilde{Q} \in \mu_{\psi^{+}}(f+P)$, using \eqref{eq:def-extended-q-neg} the next inequality holds  
\begin{equation}\label{eq:ineq-prop-2-1-previa-1}
\begin{aligned}
\int_{\{f+P>\tilde{Q}\} \cap\{Q>0\}} &\psi^{-}(|f+P-\tilde{Q}|)\,|Q| \,dx \\
&+\int_{\{f+P < \tilde{Q}\} \cap\{Q<0\}} \psi^{-}(|f+P-\tilde{Q}|)\,|Q| \,dx \\
\leqslant \int_{\{f+P < \tilde{Q}\} \cap\{Q>0\}} &\psi^{+}(|f+P-\tilde{Q}|)\,|Q| \,dx\\
&+\int_{\{f+P>\tilde{Q}\} \cap\{Q<0\}} \psi^{+}(|f+P-\tilde{Q}|)\,|Q| \,dx,
\end{aligned}
\end{equation}
for every $Q \in \Pi^{m}.$
The former inequality can be written as 
\begin{equation}\label{eq:ineq-prop-2-1}
\begin{aligned}
\int_{\{f>\tilde{Q}-P\} \cap\{Q>0\}} &\psi^{-}(|f-(\tilde{Q}-P)|)\,|Q| \,dx\\
&+\int_{\{f < \tilde{Q}-P\} \cap\{Q<0\}} \psi^{-}(|f-(\tilde{Q}-P)|)\,|Q| \,dx \\
\leqslant \int_{\{f < \tilde{Q}-P\} \cap\{Q>0\}} &\psi^{+}(|f-(\tilde{Q}-P)|)\,|Q| \,dx\\
&+\int_{\{f>\tilde{Q}-P\} \cap\{Q<0\}} \psi^{+}(|f-(\tilde{Q}-P)|)\, |Q|  \,dx,
\end{aligned}
\end{equation}
for every $Q \in \Pi^{m}.$
By \eqref{Remark} of Remark \eqref{Rem1},  \eqref{eq:ineq-prop-2-1} implies  that $\tilde{Q}-P \in \mu_{\psi^{+}}(f)$, that is, $\tilde{Q} \in$ $\mu_{\psi^{+}}(f)+P$.

Now, suppose that  $\tilde{Q}\in\mu_{\psi^{+}}(f)+P$, i.e., $\tilde{Q}-P \in \mu_{\psi^{+}}(f)$ with $f \in L^{\psi^{+}}(\Omega)$ and
therefore, inequality \eqref{eq:ineq-prop-2-1} is satisfied.
As inequality \eqref{eq:ineq-prop-2-1} is equivalent to  \eqref{eq:ineq-prop-2-1-previa-1},  then $\tilde{Q} \in$ $\mu_{\psi^{+}}(f+P)$ where $f+P \in L^{\psi^{+}}(\Omega)$.
Finally, $\mu_{\psi^{+}}(f+P)=\mu_{\psi^{+}}(f)+P$ for each $P \in \Pi^{m}$.
\end{proof}

\begin{thm}\label{thm:continuity-ebpao}
Let $\varphi \in \Phi$ such that its right derivative $\psi^+$ is strictly increasing and consider $h_{n}, h \in L^{\psi^{+}}(\Omega)$ satisfying
\begin{equation}\label{eq:hyp-continuity}
\int_{\Omega} \psi^{+}\left(\left|h_{n}-h\right|\right) \,dx \rightarrow 0
\end{equation}
as $n \rightarrow \infty$, and 
$|h_n(x)|\leqslant H(x) $ a.e. in $\Omega$  with $H \in L^{\psi^+}(\Omega).$\\
Then   $\mu_{\psi^{+}}\left(h_{n}\right) \rightarrow \mu_{\psi^{+}}(h)$ as $n \rightarrow \infty$.
\end{thm}

\begin{proof}
Let $h_{n}, h \in L^{\psi^{+}}(\Omega)$ such that \eqref{eq:hyp-continuity} holds. 

If $\mu_{\psi^{+}}\left(h_{n}\right)=P_{n}$, then by (\ref{remakb}), and using the notation of
\eqref{eq:consecuencia-caract-bpa},  we have
\[
\begin{aligned}
&\int_{\left\{h_{n}<P_{n}\right\} \cap\{Q>0\}} F^{-} (h_n, P_n, Q)
+\int_{\left\{h_{n}>P_{n}\right\} \cap\{Q<0\}} F^{-} (h_n, P_n, Q)  \\
\leqslant &\int_{\left\{h_{n}>P_{n}\right\} \cap\{Q>0\}} F^{+} (h_n, P_n, Q)
+\int_{\left\{h_{n} < P_{n}\right\} \cap\{Q<0\}} F^{+} (h_n, P_n, Q),
\end{aligned}
\]
for every $Q \in \Pi^{m}$ and for every $n \in \mathbb{N}$.

From \eqref{eq:hyp-continuity}, given $\epsilon>0$ there exists $N_{0} \in \mathbb{N}$ such that
\begin{equation}\label{eq:convergence in L1 of psi+(hn-h)}
\int_{\Omega} \psi^{+}\left(\left|h_{n}-h\right|\right) \,dx<\epsilon,
\end{equation}
for every $n>N_{0}$.
Employing the monotonicity of $\psi^{+}$, \eqref{eq:subaditiva}, \eqref{eq:convergence in L1 of psi+(hn-h)} and  the fact that $h \in L^{\psi^{+}}(\Omega)$, we get
\[
\begin{aligned}
\int_{\Omega} \psi^{+}\left(\left|h_{n}\right|\right) \,dx \leqslant \int_{\Omega} \psi^{+}\left(\left|h_{n}-h\right|+|h|\right) \,dx &\leqslant \\
\Lambda_{\psi^+}\left(\int_{\Omega} \psi^{+}\left(\left|h_{n}-h\right|\right) \,dx
+\int_{\Omega} \psi^{+}(|h|) \,dx\right)&<\Lambda_{\psi^+}
\left(\epsilon+C_{\psi^{+}, h}\right),
\end{aligned}
\]
for every $n>N_{0}$ and where $\int_{\Omega} \psi^{+}(|h|) \,dx=C_{\psi^{+}, h}<\infty$.

As $h_{n} \in L^{\psi^{+}}(\Omega)$, we have
\begin{equation}\label{eq:bound psi+ hn}
\int_{\Omega} \psi^{+}\left(|h_{n}|\right) \,dx=
C_{\psi^{+}, h_n}<\infty, \mbox { for every  } n \in \mathbb{N},
\end{equation}
we define
\[
C_{u}=\max\left\{\left\{C_{\psi^{+}, h_{n}}\right\}_{n=1}^{N_{0}}, \Lambda_{\psi^+}\left(\epsilon+C_{\psi^{+}, h}\right)\right\}<\infty,
\]
then, from \eqref{eq:bound psi+ hn} 
\[
\int_{\Omega} \psi^{+}\left(|h_{n}|\right) \,dx \leqslant C_{u}<\infty, \mbox { for every } n \in \mathbb{N}.
\]
Now we use \eqref{eq:inequality-extended} of Theorem \ref{thm:extension} to get
\begin{equation}\label{eq: bound of psi+(Pn) in L1}
\int_{\Omega} \varphi\left(\left|P_{n}\right|\right) \,dx 
\leqslant K\left\|P_{n}\right\|_{\infty} C_{u}, \mbox { for every } n \in \mathbb{N} .
\end{equation}
Next, by \eqref{eq: bound of psi+(Pn) in L1}, Jensen's inequality, inequality  (2.14) in Lemma 2.4 of \cite{CFZ12}, 
the fact that $\varphi$ is non-decreasing on $[0,\infty)$ and $\varphi \in \Delta_{2}$ with  constant $\Lambda_{\varphi}$, we obtain 
\[
\begin{aligned}
K C_{u}\left\|P_{n}\right\|_{\infty} & \geqslant \int_{\Omega} \varphi\left(\left|P_{n}\right|\right) \,dx 
\geqslant|\Omega| \varphi\left(\frac{1}{|\Omega|} \int_{\Omega}\left|P_{n}\right|\right) \,dx \\
&=|\Omega| \varphi\left(\frac{1}{|\Omega|}\left\|P_{n}\right\|_{1}\right) 
\geqslant|\Omega| \varphi\left(M_{1}\left\|P_{n}\right\|_{\infty}\right) \\
& \geqslant|\Omega| \varphi\left(\frac{1}{2^{N}}\left\|P_{n}\right\|_{\infty}\right) \geqslant \frac{|\Omega|}{\Lambda_{\varphi}^{N}} \varphi\left(\left\|P_{n}\right\|_{\infty}\right),
\end{aligned}
\]
being $M_{1}>0$ such that $M_{1}\|P\|_{\infty} \leqslant \frac{\|P\|_{1}}{|\Omega|}$ and 
$N \in \mathbb{N}$ such that $\frac{1}{2^{N}} \leqslant M_{1}$. 
Therefore,
\[
\frac{\varphi\left(\left\|P_{n}\right\|_{\infty}\right)}{\left\|P_{n}\right\|_{\infty}} \leqslant K_{u},
\]
for every $n \in \mathbb{N}$, where the constant  $K_{u}=\frac{K C_{u} \Lambda_{\varphi}^{N}}{|\Omega|}$ is independent of $P_{n}$.
\\
Then,  by \eqref{eq:des-psi+ contra phi} and the fact that 
$\psi^+(x) \rightarrow \infty$ as $x \rightarrow \infty$, 
we get that
$\left\{P_{n}\right\}$ is uniformly bounded.

Since $\left\{P_{n}\right\}$ is uniformly bounded,  then there exists a subsequence   $\left\{P_{n_{k}}\right\}$ such that $P_{n_{k}}$ converges uniformly on $\Pi^{m}$.
Let $P_{n_{k}}=\mu_{\psi^{+}}\left(h_{n_{k}}\right)$ for  $n_{k} \in \mathbb{N}$, and let   $P=\lim\limits_{n_{\mathrm{k}} \rightarrow \infty} P_{n_k}$.

By simplicity, we call  $\left\{P_{n}\right\}$ the subsequence such that  $P_{n} \rightarrow P$ as $n \rightarrow \infty$.
And, we will prove that the polynomial $P$ satisfies inequality \eqref{remakb} of Remark \ref{Rem1}.

We follow  the same way as in the proof of Theorem \ref{thm:extension}, with   $h$ and  $h_n$ 
 instead of $f$ and $f_n$,  respectively. 
 We will use the notations 
 $ F^{-} (g,S,Q)=\psi^{-} (|g-S|) |Q|,$  and
 $ F^{+} (g,S,Q)=\psi^{+} (|g-S) |Q|,$ for $g\in L^{\psi^+} (\Omega)$ and $S, Q\in \Pi^m;$
and, we will set $l=\liminf\limits_{n \rightarrow \infty}$ and $L= \limsup\limits _{n \rightarrow \infty}.$

 Thus, employing Lemma \ref{eq:conv in L1 of  psi+fn} as in the proof of Theorem \ref{thm:extension}, we have 
\[
\begin{aligned}
&\int_{\{h> P\} \cap\{Q<0\}} F^{-} (h,P,Q) +\int_{\{h<P\} \cap\{Q>0\}} F^{-} (h,P,Q)  \\
\leqslant l &\int_{\left\{h_{n}>P_{n}\right\} \cap\{Q<0\}} F^{-} (h_n,P_n,Q) 
+ l \int_{\left\{h_{n}<P_{n}\right\} \cap\{Q>0\}} F^{-} (h_n,P_n,Q)  \\
\leqslant l &\left(\int_{\left\{h_{n}>P_{n}\right\} \cap\{Q<0\}} F^{-} (h_n,P_n,Q) 
+\int_{\left\{h_{n}<P_{n}\right\} \cap\{Q>0\}} F^{-} (h_n,P_n,Q) \right) \\
\leqslant l &\left(\int_{\left\{h_{n}>P_{n}\right\} \cap\{Q>0\}} F^{+} (h_n,P_n,Q) 
+\int_{\left\{h_{n} < P_{n}\right\} \cap\{Q<0\}} F^{+} (h_n,P_n,Q) \right) \\
\leqslant L &\left(\int_{\left\{h_{n}>P_{n}\right\} \cap\{Q>0\}} F^{+} (h_n,P_n,Q)
+\int_{\left\{h_{n} < 
 P_{n}\right\} \cap\{Q<0\}} F^{+} (h_n,P_n,Q) \right) \\
\leqslant L &\int_{\left\{h_{n}>P_{n}\right\} \cap\{Q>0\}} F^{+} (h_n,P_n,Q)
+ L \int_{\left\{h_{n} < P_{n}\right\} \cap\{Q<0\}} F^{+} (h_n,P_n,Q)
\\
\leqslant  &\int_{\left\{h>P\right\} \cap\{Q>0\}} F^{+} (h,P,Q)
+  \int_{\left\{h < P\right\} \cap\{Q<0\}} F^{+} (h,P,Q),
\end{aligned}
\]
for every $Q \in \Pi^{m}$.
Thus, we have proved that $P\in \mu_{\psi^{+}} (h)$ taking into account inequality \eqref{remakb} of Remark  \ref{Rem1}. 

Finally, by Theorem \ref{thm:unique extended bpao} we have that the whole sequence $P_n$ converges to $P$.
\end{proof}

\section{Pointwise Convergence}

In this section we generalize to the Orlicz Space $L^{\psi^+}$ Corollary 2.7 of \cite{CFZ12}
and Proposition 2.6 of \cite{AFZ19}, which deal with 
the extended  best approximation operators in the setting of $L^{p}$ and $L^\varphi$, respectively.

We consider $f \in L^{\psi^+}(B)$ where $B=B(x,\varepsilon)$
is the ball in $\mathbb{R}^n,$ centered at $x$ and radius $\varepsilon>0$.
We represent the set of extended best approximation polynomials from $L^{\varphi}(B)$ to $L^{\psi^+}(B)$  by $\mu_{\psi^+}^{\varepsilon,x}(f)=\mu_{\psi^+}^B(f).$ For each polynomial in $\mu_{\psi^+}^{\varepsilon,x} (f)$  we refer to $P_x^{\varepsilon}(f).$
This polynomial can be written as
\[
P_x^{\varepsilon}(f)(t)=\sum\limits_{|\alpha|\leq m} a_{\alpha}(f)(x,\varepsilon)(t-x)^{\alpha},
\]
 being $\alpha=(\alpha_1,\dots,\alpha_n)$,\,\,  $t= (t_1, \dots, t_n),$ \,\,$t^{\alpha}=t^{\alpha_1}\dots t^{\alpha_n}$ and
 $|\alpha|=\alpha_1+\dots+\alpha_n$.

First of all, we  estimate the coefficients $a_{\alpha}(f)(x,\varepsilon)$ and
we get a pointwise convergence result for $P_x^{\varepsilon}(f)$ under a very general
differentiability assumption on the function $f$, based on the class $t_m^p(x)$ introduced
by Calder\'on and Zygmund for functions in $L^{p}$, see 
\cite[p. 172] {CZ61}.
It is worthy of mentioning that the $\Delta_2$ condition on $\psi^+$ is essential to get
the convergence mentioned above.

Before dealing with the main goal of this section, we put forward two results that will be used later.

\begin{lem}\label{cotas-cociente-delta-2}
Let $\varphi \in \Phi$ and let $\psi^+$ be its right 
derivative. 
Then, there exist positive functions $g_1,g_2$ such that
\[
g_1(\eta)\leqslant  \frac{\psi^+(\eta x)}{\psi^+(x)}\leqslant  g_2(\eta),
\]
for every $x,\eta>0$.
\end{lem}

The proof of Lemma \ref{cotas-cociente-delta-2} follows the same pattern of Lemma 2.1 in \cite{AFZ19} with $\psi^+$ instead of $\varphi$, since the continuity of $\varphi$ is not used.

Next, we get an inequality that generalizes those given in Lemma 2.4 of \cite{CFZ12}.

\begin{lem}\label{equiv-pol-phi}
Let $\varphi \in \Phi $ and let $\psi^+$ be its right 
derivative.  
Then, there exists a  constant $C_1 >0$, depending only on
$n,m$ and $\psi^+$, such that
\begin{equation}\label{des-equiv-pol-orl}
\frac{C_1}{\Lambda_{\varphi}} \psi^+\left( C_1\|P\|_{L^{\infty}_{(B(x,\varepsilon))}}\right)\leqslant 
\frac{1}{|B(x,\varepsilon)|} \int_{B(x,\varepsilon)} \psi^+(|P|)\,dt
\leqslant 
\psi^+ \left( \|P\|_{L^{\infty}_{(B(x,\varepsilon))}} \right),
\end{equation}
for all $P \in \Pi^m$.
\end{lem}

\begin{proof}
The right-hand side of inequality \eqref{des-equiv-pol-orl} follows straightforwardly.

Now, as $\varphi \in\Delta_2$, by  \eqref{eq:des-psi+ contra phi}
we get
\begin{equation}\label{delta2-p}
\frac{k}{\Lambda_{\varphi}}\|P\|_{L^{\infty}_{(B(x,\varepsilon))}}
\psi^+\left( k \|P\|_{L^{\infty}_{(B(x,\varepsilon))}}\right)
\leqslant 
\varphi\left( k\|P\|_{L^{\infty}_{(B(x,\varepsilon))}}\right),
\end{equation}
for any $k \geqslant0$. 

We use inequality (2.14) of Lemma 2.4 in \cite{CFZ12} and
the monotonicity of $\varphi$ on $[0,\infty)$ to get the existence of $C_1 \geqslant 1,$ depending only on $n$ and $m$, such that
\begin{equation}\label{equiv-monot-phi}
\varphi\left(C_1 \|P\|_{
L^{\infty}_{(B(x,\varepsilon))}}
\right)
\leqslant 
\varphi\left(\frac{1}{|B(x,\varepsilon)|} \int_{B(x,\varepsilon)}|P|\,dt\right).
\end{equation}

Using \eqref{eq:des-psi+ contra phi} and  Jensen's inequality for the convex function $\varphi$, we have 

\begin{equation}\label{jensen-p}
\begin{split}
\varphi\left(\frac{1}{|B(x,\varepsilon)|} \int_{B(x,\varepsilon)}|P|\,dt\right)
&\leqslant  \frac{1}{|B(x,\varepsilon)|} \int_{B(x,\varepsilon)}\varphi(|P|)\,dt
\\
&\leqslant  \frac{1}{|B(x,\varepsilon)|} \int_{B(x,\varepsilon)}|P|\psi^+(|P|)\,dt
\\
&\leqslant 
\|P\|_{L^{\infty}_{(B(x,\varepsilon))}}\frac{1}{|B(x,\varepsilon)|} \int_{B(x,\varepsilon)}\psi^+(|P|)\,dt.
\end{split}
\end{equation}

Then, from \eqref{delta2-p}, \eqref{equiv-monot-phi} with $k=C_1$ and  \eqref{jensen-p},  we get the required inequality. 
\end{proof}

Recall that if  $P_x^{\varepsilon}(f)\in \mu_{\psi^+}^{\varepsilon,x}(f)$, 
by Theorem \ref{thm:extension},  we have 
\begin{equation}\label{fuerte-1-o}
\int_{B(x,\varepsilon)} \varphi(|P_x^{\varepsilon}|)\,dt \leqslant 
5\Lambda_{\psi^+} \|P_x^{\varepsilon}\|_{L^{\infty}_{(B(x,\varepsilon))}} \int_{B(x,\varepsilon)} \psi^+(|f|)\,dt.
\end{equation}

Then, note that,  if we consider \eqref{delta2-p} with $k=C_1,$ inequality  \eqref{equiv-monot-phi} and Jensen's Inequality, we obtain
\[
\frac{C_1}{\Lambda_{\varphi}}\|P_{x}^{\varepsilon} \|_{L^{\infty}_{(B(x,\varepsilon))}}
\psi^+\left( C_1 \|P_{x}^{\varepsilon}\|_{L^{\infty}_{(B(x,\varepsilon))}}\right)
\leqslant \frac{1}{|B(x,\varepsilon)|} \int_{B(x,\varepsilon)} \varphi\left(|P_{x}^{\varepsilon}|\right)\,dt,
 \]
and from  \eqref{fuerte-1-o} we get
\begin{equation}\label{P-eps-inf-prom}
\psi^+\left(C_1\|P_x^{\varepsilon}\|_{L^{\infty}_{(B(x,\varepsilon))}}\right)\leqslant 
\frac{5\Lambda_{\psi^+} \Lambda_{\varphi}}{C_1} \frac{1}{|B(x,\varepsilon)|}  \int_{B(x,\varepsilon)} \psi^+(|f|)\,dt,
\end{equation}
for $f \in L^{\psi^+}(B(x,\varepsilon))$ and where $C_1$ is a constant independent of $x,\varepsilon$  and $f$.

We point out that inequality \eqref{P-eps-inf-prom} is analogous 
to (2.21) in \cite{CFZ12}.

Until this point, we have estimated the extended best approximation 
operator. Next, we estimate the coefficients of $P_x^{\varepsilon}(f),$ for $f \in L^{\psi^+}_{loc}(\mathbb{R}^n),$ where $\psi^+$ is not necessarily a continuous function, as it has done in Corollary 2.5 in \cite{CFZ12} and Proposition 2.3 of \cite{AFZ19}, for the instances of $L^{p-1}$ with $p>1$ and $L^\varphi,$ where $\varphi$ was a continuous and non-decreasing function, respectively.

\begin{prop}\label{cota-coef-prom}
Let $\varphi \in \Phi $ and let $\psi^+$ be its right 
derivative. 
Besides, suppose that $f \in L^{\psi^+}_{loc}(\mathbb{R}^n)$. Then
\[
\psi^+ \left( \frac{|\varepsilon|^{|\alpha|}}{C}  |a_{\alpha}(f)(x,\varepsilon)|\right)
\leqslant 
\frac{K}{|B(x,\varepsilon)|} \int_{B(x,\varepsilon)} \psi^+(|f|)\,dt,
\]
for every $\varepsilon>0$, $|\alpha| \leqslant  m$ and $x \in \mathbb{R}^n$ and where the constants $C$ and $K$ depend only on $n,m$ and $\psi^+$.
\end{prop}

\begin{proof}
For each  $P_x^{\varepsilon}(f) \in \mu_{\psi^+}^{\varepsilon,x}(f)$, from Lemma 2.4 in \cite{CFZ12}, there exists $C>0$ depending on $n$ and $m$ such that
\[
\frac{1}{C}
\|P_x^{\varepsilon}\|_{L^{\infty}_{(B(x,\varepsilon))}} 
\leqslant 
\max\limits_{|\alpha|\leqslant  m} 
\varepsilon ^{|\alpha|} |a_{\alpha}(f)(x,\varepsilon)|
\leqslant  C \|P_x^{\varepsilon}\|_{L^{\infty}_{(B(x,\varepsilon))}}.
\]
As a consequence of \eqref{P-eps-inf-prom} and the fact that $\psi^+$ is a non-decreasing function, we obtain
\[
\begin{split}
\psi^+\left( \frac{1}{C} \max\limits_{|\alpha|\leq m}
\varepsilon ^{|\alpha|}|a_{\alpha}(f)(x,\varepsilon)|  \right) 
&\leqslant 
\psi^+\left(\|P_x^{\varepsilon}\|_{L^{\infty}_{(B(x,\varepsilon))}}\right)
\\
\leqslant 
\Lambda_{\psi^{+}} \psi^+\left(C_1 \|P_x^{\varepsilon}\|_{L^{\infty}_{(B(x,\varepsilon))}}\right)
&\leqslant
\frac{K}{|B(x,\varepsilon)|}\int_{B(x,\varepsilon)}\psi^+(|f|)\,dt,
\end{split}
\]
for a suitable constant $K>0.$
\end{proof}

Now we define a pointwise smoothness condition, according to that introduced by Calderón and Zygmund in \cite{CZ61}.

\begin{defi}
Let  $\varphi \in \Phi$ and let $\psi^+$ be its right 
derivative. Assume that $f \in L^{\psi^+}_{loc}(\mathbb{R}^n)$.
\\
We say that $f \in t_m^{\psi^+}(x)$ if there exists a polynomial $P_x(f)(t)$ of degree at most $m$
such that
\begin{equation}\label{condicion-t-phi-m}
\frac{1}{|B(x,\varepsilon)|} \int_{B(x,\varepsilon)}\psi^+(|f(t)-P_x(f)(t)|)\,dt=
o(\psi^+(\varepsilon^m)) \;\;\mbox{as}\;\varepsilon \to 0.
\end{equation}
\end{defi}

The uniqueness of the polynomial $P_x(f)$ follows from 
condition \eqref{condicion-t-phi-m} as we shall see next. 

\begin{thm}\label{thm:unicidad-map-extendido-bola}
Let $\varphi \in \Phi$ and let $\psi^+$ be its right 
derivative. Then, for $f \in t_m^{\psi^+}(x),$ there exists a unique polynomial $P_x(f)$ that satisfies \eqref{condicion-t-phi-m}.
\end{thm}

\begin{proof}
Assume that $P^1_x(f)$ and $P^2_x(f)$ are
polynomials that satisfy  \eqref{condicion-t-phi-m}.
Then, by \eqref{eq:subaditiva}, we get
\begin{equation}\label{condicion-t-phi-m-diferencia}
\frac{1}{|B(x,\varepsilon)|} 
\int_{B(x,\varepsilon)}\frac{\psi^+(|P_x^1(f)(t)-P_x^2(f)(t)|)}{\psi^+(\varepsilon^m)}\,dt\to 0
\;\;\mbox{as}\;\varepsilon \to 0.
\end{equation}

Let $Q(f)(t)=P_x^1(f)(t)-P_x^2(f)(t)$.
We will see that such a polynomial  $Q(f) \in\Pi^m$ that satisfies \eqref{condicion-t-phi-m-diferencia}
has to be the zero polynomial.

As $\psi^+ \in \mathcal{S}$, by \eqref{condicion-t-phi-m-diferencia},
Lemma \ref{equiv-pol-phi} and Lemma 2.4 iii) in \cite{CFZ12},
there exists $C>0$, depending only on $n,m$ and $\psi^+$, such that
\[
\frac{\psi^+ \left(\varepsilon^{|\alpha|}\frac{|b_{\alpha}|}{C}\right)}{\psi^+(\varepsilon^{|\alpha|})}
\leqslant 
\frac{\psi^+ \left(\varepsilon^{|\alpha|}\frac{|b_{\alpha}|}{C}\right)}{\psi^+(\varepsilon^m)}
\leqslant  \frac{\psi^+\left(\|Q(f)\|_{L^{\infty}_{(B(x,\varepsilon))}}\right)}{\psi^+(\varepsilon^m)}
\to 0
\;\;\mbox{as}\;\varepsilon \to 0,
\]
for every $\alpha$ such that $|\alpha|\leqslant  m$ and where
$Q(f)(t)=\sum\limits_{|\alpha|\leq m}b_{\alpha}(x)(t-x)^{\alpha}$.
Now, by Lemma \ref{cotas-cociente-delta-2} and the fact that $\varphi \in \Phi$, 
we have $b_{\alpha}=0$ for every  $\alpha$ such that $|\alpha|\leqslant  m$,
that is, $Q(f)=P_x^1(f)-P_x^2(f)\equiv 0$ on $B(x,\varepsilon)$.
\end{proof}

From now on, we  write 
\[
P_x(f)(t)=\sum\limits_{|\alpha|\leq m} a_{\alpha}(f)(x) (t-x)^{\alpha},
\]
and we denote by $\partial ^{\alpha}f(x)$ the coefficients of the polynomial $P_x(f)$ multiplied by $\alpha!$.

Next, we get a pointwise convergence result for the extended best $\varphi$-approximation polynomial.

\begin{prop}\label{conv-puntual}
Let $\varphi \in \Phi $ with right derivative function $\psi^+$ and  let $f \in L^{\psi^+}_{loc}(\mathbb{R}^n)$.
Additionally, suppose that $\psi^+$ is a strictly increasing function and
$f \in t_m^{\psi^+}(x)$ in a fixed point $x$.
\\
Then $P_x^{\varepsilon}(f)(t)  \to P_x(f)(t)$ as $\varepsilon \to 0$ for every $t\in\mathbb{R}^n$.
\end{prop}

\begin{proof}
As $\psi^+$ is a strictly increasing function, 
then $\varphi(x)=\int_0^x \psi^+(t)\,dt$ is also  strictly convex. 
From Theorem \ref{thm:prop-traslacion}, we have
$P_x^{\varepsilon}(f-P_x(f))=P_x^{\varepsilon}-P_x(f)$. 
As $\psi^+$ is a strictly increasing function, by Theorem \ref{thm:unique extended bpao}, 
$P_x^{\varepsilon}(f-P_x(f))$ and $P_x^{\varepsilon}(f)$
are the unique polynomials  of the sets $\mu_{\psi^+}(f-P_x(f))$ and $\mu_{\psi^+}(f)$ respectively, then
$a_{\alpha}(f-P_x(f))(x,\varepsilon)=a_{\alpha}(f)(x,\varepsilon)-\frac{\partial ^{\alpha}f(x)}{\alpha!}$.

By Proposition \ref{cota-coef-prom}, we get
\[
\psi^+\left( \frac{1}{C} \varepsilon^{|\alpha|} \left|a_{\alpha}(f)(x,\varepsilon)-\frac{\partial ^{\alpha}f(x)}{\alpha!}\right|\right)\leqslant 
\frac{K}{|B(x,\varepsilon)|}\int_{B(x,\varepsilon)} \psi^+(|f(t)-P_x(f)(t)|)\,dt.
\]
And $f \in t_m^{\psi^+}(x)$, then
\[
\frac{
\frac{1}{|B(x,\varepsilon)|}\int_{B(x,\varepsilon)} \psi^+(|f(t)-P_x(f)(t)|)\,dt}{\psi^+(\varepsilon^m)}\to 0\;\;
\mbox{as}\;\varepsilon \to 0.
\]
Therefore, 
\begin{equation}\label{limite-total}
\frac{\psi^+\left( \frac{1}{C} \varepsilon^{|\alpha|} \left|a_{\alpha}(f)(x,\varepsilon)-\frac{\partial ^{\alpha}f(x)}{\alpha!}\right|\right)}{\psi^+(\varepsilon^m)} \to 0 \;\;\mbox{as}\;\;\varepsilon \to 0, \;\mbox{for}\; |\alpha|\leqslant  m.
\end{equation}
We shall see that \eqref{limite-total} implies that $\left|a_{\alpha}(f)(x,\varepsilon)-\frac{\partial ^{\alpha}f(x)}{\alpha!}\right|\to 0$ as $\varepsilon \to 0$. To do so, we write
\[
\frac{\psi^+\left( \frac{1}{C} \varepsilon^{|\alpha|}
\left|a_{\alpha}(f)(x,\varepsilon)-\frac{\partial ^{\alpha}f(x)}{\alpha!}\right|\right)}{\psi^+(\varepsilon^m)} = \]
\[
\frac{\psi^+\left( \frac{1}{C} \varepsilon^{|\alpha|}
\left|a_{\alpha}(f)(x,\varepsilon)-\frac{\partial ^{\alpha}f(x)}{\alpha!}\right|\right)}{\psi^+(\varepsilon^{|\alpha|})}\frac{\psi^+(\varepsilon^{|\alpha|})}{\psi^+(\varepsilon^m)};
\]
and, due to $1\leqslant \frac{\psi^+(\varepsilon^{|\alpha|})}{\psi^+(\varepsilon^m)}$ as $\varepsilon \to  0$ for $|\alpha| \leqslant  m$, then
\begin{equation}\label{limite-cociente}
\frac{\psi^+\left( \frac{1}{C} \varepsilon^{|\alpha|}
\left|a_{\alpha}(f)(x,\varepsilon)-\frac{\partial ^{\alpha}f(x)}{\alpha!}\right|\right)}
{\psi^+(\varepsilon^{|\alpha|})}\to 0 \;\;\mbox{as}\;\varepsilon \to  0\;\;
\mbox{and for}\;\;|\alpha|\leqslant  m.
\end{equation}

Suppose that $\frac{1}{C} \left|a_{\alpha}(f)(x,\varepsilon)-\frac{\partial ^{\alpha}f(x)}{\alpha!}\right|\not\to 0 $.
Then, there exist $\eta>0$  and a sequence $\{\varepsilon_l\}$ such that $\varepsilon_l\to 0$ and
$\frac{1}{C} \left|a_{\alpha}(f)(x,\varepsilon_l)-\frac{\partial ^{\alpha}f(x)}{\alpha!}\right|\geqslant  \eta>0$.
Since $\psi^+$ is a strictly increasing function, we have
\[
\frac{\psi^+\left( \frac{1}{C} \varepsilon_l^{|\alpha|} \left|a_{\alpha}(f)(x,\varepsilon_l)
-\frac{\partial ^{\alpha}f(x)}{\alpha!}\right|\right)}{\psi^+\left(\varepsilon_l^{|\alpha|}\right)}>
\frac{\psi^+\left(\varepsilon_l^{|\alpha|} \eta \right)}{\psi^+\left({\varepsilon_l^{|\alpha|}}\right)},
\]
for every $l>0$ and for $|\alpha|\leq m$.

As $\psi^+$ satisfies the $\Delta_2$ condition, Lemma \ref{cotas-cociente-delta-2} holds; consequently,
there exists a positive function $g_1$
such that
\[
\frac{\psi^+\left(\varepsilon_l^{|\alpha|} \eta \right)}
{\psi^+\left({\varepsilon_l^{|\alpha|}}\right)} \geqslant 
g_1(\eta)>0  \;\;\mbox{for every}\;\; \varepsilon_l>0 \,\; \mbox{and for  }\;\; |\alpha|\leq m,
\]
then
\begin{equation}\label{cota-inferior-cociente}
\frac{\psi^+\left(\varepsilon_l^{|\alpha|} \eta \right)}
{\psi^+\left({\varepsilon_l^{|\alpha|}}\right)} \not \to 0  \;\;\mbox{as}\;\; \varepsilon_l \to 0,\;\;\mbox{for}\;\;|\alpha|\leq m.
\end{equation}

Now, from  \eqref{limite-cociente} and \eqref{cota-inferior-cociente}, we reach a contradiction.  
Thus,
\begin{equation}\label{conv-coef}
\left|a_{\alpha}(f)(x,\varepsilon)-\frac{\partial ^{\alpha}f(x)}{\alpha!}\right|\to 0\;\; \mbox{as} \;\;\varepsilon \to 0
\,\; \mbox{for  }\;\; |\alpha|\leq m.
\end{equation}

Last, \eqref{conv-coef} implies that $P_x^{\varepsilon}(f)(t)  \to P_x(f)(t)$ for each $t\in \mathbb{R}^n$ as $\varepsilon \to 0$.
\end{proof}

In the event where $\psi^+$ is a strictly increasing function, we have just proved the pointwise convergence of the extended best polynomial approximation for functions of $L^{\psi^+}(\mathbb{R}^n)$ to the polynomial that satisfies \eqref{condicion-t-phi-m}. This outcome is a generalization of Proposition 2.6 in \cite{AFZ19}, which makes the assumption of differentiability on $\varphi.$ The convergence of the best approximation for the situation $\varphi (x) = x^p,$ was obtained in Corollary 2.7 in \cite{CFZ12}, which applies the traditional differentiability criterion that was first proposed by Calderón and Zygmund in \cite{CZ61}.

\section{Acknowledgements}
We value the referees' comprehensive and perceptive feedback, which improved this work.

\end{document}